\documentclass[11pt,a4paper]{article}
\usepackage[francais,english]{babel}
\usepackage{amsmath}
\usepackage{amssymb}
\usepackage{latexsym}
\usepackage{amsfonts}

\newcommand{\NN}{{\mathbb N}}

\newcommand{\CC}{{\mathbb C}}
\newcommand{\ZZ}{{\mathbb Z}}
\newcommand{\bs}{\bigskip}
\newcommand{\q}[1]{\left[ #1\right] }

\renewcommand{\deg}[1]{{\partial}^{o}#1}

\renewcommand{\a}{\alpha}
\newcommand{\ab}{{\a}_1}
\newcommand{\ac}{{\a}_2}
\newcommand{\ai}{{\a}_{2k-2}}
\newcommand{\aj}{{\a}_{2k-1}}
\newcommand{\af}{{\a}_{2N-2}}
\newcommand{\ag}{{\a}_{2N-1}}
\renewcommand{\b}{\beta}
\newcommand{\bc}{{\b}_2}
\newcommand{\bi}{{\b}_{2k-2}}
\newcommand{\bj}{{\b}_{2k-1}}

\newcommand{\bg}{{\b}_{2N-1}}
\newcommand{\ob}{(x_1 y_1)}
\newcommand{\oc}{(y_1 x_2)}
\newcommand{\oi}{(y_{k-1} x_k)}
\newcommand{\oj}{(x_k y_k)}
\newcommand{\of}{(y_{N-1} x_{N})}

\newcommand{\hbab}{{\ob}^{-(\ab}-1)}
\newcommand{\hb}{{\ob}^{-\ab}}
\newcommand{\hc}{{\oc}^{-\ac}}
\newcommand{\hi}{{\oi}^{-\ai}}
\newcommand{\hj}{{\oj}^{-\aj}}
\newcommand{\hf}{{\of}^{-\af}}

\newcommand{\gc}{{\oc}^{-\bc}}
\newcommand{\gi}{{\oi}^{-\bi}}
\newcommand{\gj}{{\oj}^{-\bj}}

\newcommand{\HT}{\hb\cdots\hf}

\newcommand{\HV}{{{\ob}^{-1}}\hc\ldots\hf}
\newcommand{\C}[2]{\left[ \begin{array}{c}#1\\#2\end{array}\right] }
\newcommand{\cnp}[2]{\left( \begin{array}{c}#1\\#2\end{array}\right) }

\renewcommand{\dfrac}[2]{\displaystyle{\frac{#1}{#2}}}
\newcommand{\tab}{\hspace*{\fill}}
\newcommand{\ok}{\tab $\Box$\bs}

\newcommand{\SIA}{\displaystyle{
\sum\limits_{{\underline{\a}}\in J_n}}}
\newcommand{\SIB}{\displaystyle{
\sum\limits_{{\underline{\b}}\in J_n}}}
\newcommand{\SIC}{\displaystyle{
\sum\limits_{{\underline{\a}}\in J'_{n-1}}}}
\newcommand{\SID}{\displaystyle{
\sum\limits_{{\underline{\a}}\,\in J''_{n+1}\atop{\ab =2}}}}
\newcommand{\SJ}{\displaystyle{
\sum\limits_{{\underline{\a}}\in J_n\atop{\ab =1}}}}
\newcommand{\SK}{\displaystyle{
\sum\limits_{{\underline{\a}}\in J_n\atop{\ab >1}}}}
\newcommand{\SP}{{\Sigma}^{+}}
\newcommand{\SM}{{\Sigma}^{-}}
\newtheorem{lem}{Lemme}[section]

\newtheorem{thm}{Th{\'e}or{\`e}me}[section]
\newtheorem{cor}{Corollaire}[section]

\newenvironment{Def}{{\noindent\bf D{\'e}finitions. }}{}
\newenvironment{ddf}{{\noindent\bf D{\'e}finition. }}{}
\newenvironment{dem}{{\noindent\bf D{\'e}monstration. }}{\ok}
\newenvironment{pr}{{\noindent\bf Preuve. }}{}

\title{Sur les int{\'e}grales de mouvement du syst{\`e}me de
sinus-Gordon sur r{\'e}seau}
\author{C.Grunspan\\
Ecole Polytechnique\\
Centre de math{\'e}matiques\\
UMR 7640 du CNRS\\
F-91128 PALAISEAU Cedex\\
grunspan\at math.polytechnique.fr}
\date{}
\numberwithin{equation}{section}

\begin{document}
\maketitle

\selectlanguage{french}
\begin{abstract}
\noindent
Nous donnons des formules explicites pour des 
densit{\'e}s d'int{\'e}grales de mouvement du syt{\`e}me 
de sinus-Gordon sur r{\'e}seau (quantique ou non). 
La s{\'e}rie g{\'e}n{\'e}ratrice de ces densit{\'e}s d'int{\'e}grales 
de mouvement peut {\^e}tre vue comme le d{\'e}veloppement
du logarithme d'une fraction continue ({\'e}ventuellement quantique).
Dans le cas o{\`u} $q$ est une racine de l'unit{\'e}, on montre que 
certaines int{\'e}grales de mouvement s'identifient aux int{\'e}grales de
mouvement classiques.
\end{abstract}

\bs
\selectlanguage{english}
\begin{abstract}
\noindent
We give explicit formulas for some densities
of integrals of motion for the sine-Gordon system on lattice
(quantum or not). The generating function for the densities of
integrals of motion may be seen as the expansion of the logarithm of a
certain continuous fraction (possibly quantum).
In the case of $q$ root of the unity, we show that these integrals of
motion can be identified to the classical integrals of motion.
\end{abstract}

\bs
\selectlanguage{french}
\quad
{\bf Mots Cl{\'e}s~:} Th{\'e}orie de Toda sur r{\'e}seau, coefficients
$q$ binomi-\hfill

\quad aux.

\selectlanguage{english}
\quad
{\bf Key words~:} Toda theory on lattice, $q$ binomial coefficients.

\newpage
\selectlanguage{french}
\setcounter{section}{-1}
\section{Introduction}
\noindent
Le probl{\`e}me des int{\'e}grales de mouvement de la th{\'e}orie de Toda 
continue a {\'e}t{\'e} introduit par A. Zamolodchikov en liaison avec 
des probl{\`e}mes de th{\'e}orie int{\'e}grale des champs et de th{\'e}orie 
conforme (\cite{ZALZA}, \cite{ZA}, \cite{ZAM}).
Il a ensuite {\'e}t{\'e} {\'e}tudi{\'e} par R. Sasaki et I. Yamanaka sous
l'angle de la th{\'e}orie de KdV (\cite{SAS}). Il consiste {\`a} trouver 
des op{\'e}rateurs qui commutent {\`a} des op{\'e}rateurs d'{\'e}cran
dans une ${\cal W}-$alg{\`e}bre classique
ou une Alg{\`e}bre d'Op{\'e}rateurs de Vertex, dans le cas quantique.
Ce probl{\`e}me a {\'e}t{\'e} r{\'e}solu par B. Feigin et E. Frenkel par des
m{\'e}thodes cohomologiques (\cite{FEI}). Malheureusement, il n'existe pas
v{\'e}ritablement de formules explicites pour ces int{\'e}grales de
mouvement.

D'autre part, une version sur r{\'e}seau de la th{\'e}orie de Toda 
a {\'e}t{\'e}
introduite par A.G. Izergin et V.E. Korepin (\cite{IZK}, \cite{IZKO}),
et depuis {\'e}tudi{\'e}e par
de nombreux auteurs (voir par exemple \cite{BABO}, \cite{BOB},
\cite{FAD}, \cite{RES}).
Au niveau classique, B. Enriquez et B. Feigin ont montr{\'e} {\`a} l'aide
de d{\'e}singularisation de Demazure et de m{\'e}thodes cohomologiques,
que l'on pouvait obtenir des int{\'e}grales de mouvement pour le mod{\`e}le
de sinus-Gordon, en
d{\'e}veloppant le logarithme d'une 
certaine fraction continue (\cite{ENR} et th{\'e}or{\`e}me 3.1 ci-dessous).
Nous nous proposons de g{\'e}n{\'e}raliser ce r{\'e}sultat au cas
quantique, de montrer que les int{\'e}grales de mouvement obtenues sont
quantifiables, et de donner des formules explicites.
Nous appliquerons ces formules dans le cas o{\`u} $q$ est une racine de
l'unit{\'e}, et nous verrons que l'on retrouve les formules classiques.
Ce ph{\'e}nom{\`e}ne est {\`a} rapprocher de l'homomorphisme de
Frobenius-Lusztig.

Un probl{\`e}me ouvert par notre travail serait de prouver la
commutativit{\'e} des int{\'e}grales de mouvement (comme
dans le cas continu), et de calculer le spectre de ces
op{\'e}rateurs aux racines de l'unit{\'e}.

\section{Pr{\'e}sentation du syst{\`e}me de sinus-Gordon 
quantique sur r{\'e}seau}

\noindent
Soient $n\in{\NN}^*$ un entier fix{\'e}, et $q$ une ind{\'e}termin{\'e}e. 
On s'int{\'e}resse {\`a} l'alg{\`e}bre $A_q$ engendr{\'e}e sur  
$\CC[q,q^{-1}]$ 
par les $x_i^{\pm 1}$ et $y_i^{\pm 1}, i\in\ZZ$
soumis aux relations~:
$x_ix_j=qx_j x_i,\, y_iy_j=qy_jy_i$,
$x_iy_j=q^{-1}y_jx_i,\, y_ix_j=q^{-1}x_jy_i$ pour $i<j$,
et $x_iy_i=q^{-1}y_ix_i$.

En utilisant par exemple des extensions de Ore (\cite{KAS}),
on montre que $A_q$ est un anneau int{\`e}gre,
et que c'est un $\CC[q,q^{-1}]$ module libre
dont une base est donn{\'e}e par la famille 
$\prod\limits_{-\infty}^{+\infty}x_i^{\a_i}y_i^{\b_i}$, 
o{\`u} $(\a_i)$ et $(\b_i)$ sont deux suites presque nulles dans
${\ZZ}^{\ZZ}$.

L'alg{\`e}bre $A_q$ est gradu{\'e}e en posant~:
    $$
      \forall\, i\in\ZZ,\quad \deg{x_i}=-\deg{y_i}=1.
    $$
On note $A_q[0]$ la partie de $A_q$ de degr{\'e} $0$.
C'est la sous-alg{\`e}bre de $A_q$ engendr{\'e}e par les "{\'e}pingles''
$(x_i y_i)$ et $(y_i x_{i+1}),\, i\in\ZZ$.

Si $a$ et $b$ dans $A_q$ sont homog{\`e}nes, on pose~:
    $$
      [a,b]_q=ab-q^{(\deg{a})(\deg{b})}ba.
    $$
De mani{\`e}re g{\'e}n{\'e}rale, si $a=\sum_n a_n$ et
$b=\sum_p b_p$, avec $a_n$ et $b_p$ homog{\`e}nes, on pose~:
    $$
      [a,b]_q=\sum_{n,p}\, [a_n,b_p]_q.
    $$
Ce $q-$crochet d{\'e}finit {\`a} la limite pour $q\rightarrow 1$, une
structure de Poisson sur $\CC[x_i,y_i,\, i\in\ZZ]$ par~:
\begin{align}
\forall i\in\ZZ,\quad \{ x_i,y_i\} &=-x_i y_i,\\
\forall i<j,\quad \{ x_i,x_j\} &=x_i x_j,\\
\{ x_i,y_j\} &=-x_i y_j,\\
\{ y_i,y_j\} &=y_i y_j,\\
\{ y_i,x_j\} &=-y_i x_j.
\end{align}
Les op{\'e}rateurs d'{\'e}cran $\Sigma^{+}=\displaystyle
\sum\limits_{i=-\infty}^{+\infty} x_i$ et 
$\Sigma^{-}=\displaystyle
\sum\limits_{i=-\infty}^{+\infty} y_i$ n'appartiennent pas {\`a}
$A_q$.
Cependant, les morphismes~:
   $$
    \begin{array}{rcl}
    V_{+}\quad: A_q[0]&\longrightarrow&A_q\\
      P&\longmapsto&[\Sigma^{+},P]=
      \sum\limits_{i=-\infty}^{+\infty} [x_i,P]
    \end{array}
   $$
et
   $$
    \begin{array}{rcl}
    V_{-}\quad: A_q[0]&\longrightarrow&A_q\\
      P&\longmapsto&[\Sigma^{-},P]=
      \sum\limits_{i=-\infty}^{+\infty} [y_i,P]
    \end{array}
   $$
sont bien d{\'e}finis car $x_i$ (resp. $y_i$) commute avec toute
\'epingle ne contenant pas de $x_i$ (resp. $y_i$).

\section{Int{\'e}grales de mouvement pour le syst{\`e}me 
de sinus-Gordon sur r{\'e}seau}

\noindent 
Commen\c cons par rappeler la d{\'e}finition d'une densit{\'e} d'int{\'e}grale 
de mouvement.

Soit $T$ l'application de translation sur $A_q$~:
   $$
    \begin{array}{rcl}
    T:\quad A_q&\longrightarrow&A_q\\
             x_i&\longmapsto&x_{i+1}\\
             y_i&\longmapsto&y_{i+1}
    \end{array}
   $$
Il est clair que $T$ est un automorphisme d'alg{\`e}bre gradu{\'e}.

\bs
  \begin{ddf}
   On dit que $P\in A_q[0]$ est une densit{\'e} d'int{\'e}grale de 
   mouvement du syst{\`e}me
   de sinus-Gordon quantique sur r{\'e}seau 
   si $V_{+}.P$ et $V_{-}.P$ appartiennent {\`a} 
   $\textnormal{Im}(T-Id)$.
  \end{ddf}
 
\bs
Etant donn{\'e} une densit{\'e} d'int{\'e}grale de mouvement $P$,
on conviendra que l'int{\'e}grale de mouvement associ{\'e}e est 
$\sum_{k\in\ZZ} T^k P$. C'est un {\'e}l{\'e}ment d'une certaine
compl{\'e}tion de $A_q$ qui commute avec les op{\'e}rateurs
d'{\'e}cran $\SP$ et $\SM$.

Nous allons donner des formules explicites pour des densit{\'e}s 
d'int{\'e}grales de mouvement que nous noterons $\psi_n$.
Pour ce faire, il nous faut d{\'e}finir les $q-$nombres et
factorielles quantiques.

\bs
 \begin{Def}
   On d{\'e}finit successivement~:
  \begin{itemize}
   \item pour $n\in\NN,\, [n]=\displaystyle{q^n -1\over q-1}
        =\sum_{k=0}^{n-1} q^k\in\ZZ[q]$\,;
   \item pour $n\in\NN,\, [n]!=\displaystyle\prod\limits_{k=1}^{n}\, [k]$\,;
   \item pour $(n,p)\in{\ZZ}^2$,
        $$
         \C{n}{p}=\left\{
          \begin{array}{ll}
           1,&\text{si }p=0\,;\\
           0,&\text{si }p<0\text{ ou si }p>\text{Max}(0,n)\,;\\
           \displaystyle{[n]!\over [p]! [n-p]!},&\text{si }0\leq p\leq n\,;
          \end{array}
         \right.
        $$
   \item pour des entiers $N,a_1,\ldots,a_N$, 
        $$
        F_q(a_1,\ldots,a_N)=\displaystyle\prod\limits_{i=1}^{N-1}
        \C{a_i+a_{i+1}-1}{a_{i+1}}.
        $$
  \end{itemize}

 \end{Def}

\bs
Avant de d{\'e}finir les $\psi_n$, nous avons besoin de 
d{\'e}montrer le lemme suivant~:

   \begin{lem}\label{lemint}
    Soient $n,\, a_1,\ldots,a_n$ des entiers avec $a_1>0$. Alors, 
     $$
      \displaystyle{\frac{{[}a_1+\ldots +a_n{]}}{{[}a_1{]}}}
       F_q(a_1,\ldots,a_n)\in \ZZ[q].
     $$
   \end{lem}

\begin{dem}
La d{\'e}monstration se fait par r{\'e}currence sur $n$.
Si $n=1$, il n'y a rien {\`a} d{\'e}montrer, et si $n=2$, alors~: 
  $$
   \displaystyle{\frac{{[}a_1+
    a_2{]}}{{[}a_1{]}}}F_q(a_1,a_2)=\C{a_1+a_2}{a_1}\in\ZZ[q].
  $$
Supposons $n>2$. 

Si $a_2=0$, alors $F_q(a_1,\ldots,a_n)=0$ sauf si
$$
a_2=\ldots =a_n=0,
$$ 
auquel cas, $F_q(a_1,\ldots,a_n)=1$ et le lemme est
clair. 

Si $a_2>0$, on a les identit{\'e}s~:
\begin{eqnarray*}
{[}a_1+a_2\ldots +a_n{]}&=&q^{a_1}{[}a_2+\ldots +a_n{]}+{[}a_1{]}\\
F_q(a_1,\ldots,a_n)&=&\C{a_1+a_2-1}{a_2}F_q(a_2,\ldots,a_n)\\
\displaystyle{\frac{1}{{[}a_1{]}}}\C{a_1+a_2-1}{a_2}&=&
\displaystyle{\frac{1}{{[}a_2{]}}}\C{a_1+a_2-1}{a_1}
\end{eqnarray*}
D'o{\`u} l'on d{\'e}duit que 
\begin{multline*}
\displaystyle{\frac{{[}a_1+\ldots +a_n{]}}{{[}a_1{]}}}
F_q(a_1,\ldots,a_n)\\
=F_q(a_1,\ldots,a_n)+q^{a_1}\C{a_1+a_2-1}{a_1}
\displaystyle{\frac{{[}a_2+\ldots +a_n{]}}{{[}a_2{]}}}
F_q(a_2,\ldots, a_n).
\end{multline*}
On applique ensuite l'hypoth{\`e}se de r{\'e}currence {\`a}
$a_2,\ldots,a_n$ qui forment une suite de $n-1$ entiers
avec $a_2>0$. D'o{\`u} le r{\'e}sultat. 
\end{dem}

Nous en avons presque fini avec les d{\'e}finitions et lemmes auxiliaires.
Notons encore $T^{{1\over 2}}$ l'application de demi-translation~:
   $$
    \begin{array}{rcl}
    T^{{1\over 2}}:\quad A_q&\longrightarrow&A_q\\
                         x_i&\longmapsto&y_i\\
                         y_i&\longmapsto&x_{i+1}
    \end{array}
   $$
C'est un automorphisme d'alg{\`e}bre (non gradu{\'e}).

Nous sommes maintenant en mesure de d{\'e}finir les $\psi_n$.

\bs
\noindent
{\bf D{\'e}finition des $\psi_n$.}
Pour $n\in{\NN}^{*}$, on pose~:
   \begin{align}
     \psi_n&:=A_n+B_n,\\
     \intertext{avec:}
     B_n&:=T^{{1\over 2}}A_n,
   \end{align}
et~:
    \begin{equation}
      A_n:=\displaystyle\sum\limits_{{\underline\a}\in I_n}
      \displaystyle{[n]\over [\alpha_1]}
      F_q(\alpha_1,\ldots,\alpha_{2N-2})
      (x_1 y_1)^{-\alpha_1}\ldots (y_{N-1}x_N)^{-\alpha_{2N-2}},
    \end{equation}
o{\`u} $N$ est un entier quelconque v{\'e}rifiant~: 
    \begin{equation}\label{ineq}
      n\leq 2(N-1),
    \end{equation}
et $I_n$ est l'ensemble des multi-indices 
    \begin{equation}
     {\underline\a}=(\alpha_1,\ldots,\alpha_{2N-2})\in{\NN}^{*}\times
     {\NN}^{2N-3}
    \end{equation}
v{\'e}rifiant~:
    \begin{equation}
      \alpha_1+\ldots+\alpha_{2N-2}=n.
    \end{equation}

\bs
D'apr{\`e}s le lemme pr{\'e}c{\'e}dent, $\psi_n$ est bien d{\'e}finie. De plus, on
peut voir que $\psi_n$ ne d{\'e}pend pas de $N$ v{\'e}rifiant l'in{\'e}galit{\'e}
\ref{ineq}.
Ceci provient du fait
suivant d{\'e}ja utilis{\'e} dans la d{\'e}monstration du lemme \ref{lemint}~:

\bs
\noindent
{\bf Fait 1~:} Si $m,\a_1,\ldots,\a_i,\ldots,\a_m$ sont des entiers
tels que $\a_i=0$, alors~:
   $$
    F_q(\a_1,\ldots,\a_m)\not=0\,\Longrightarrow\, 
    \a_i=\a_{i+1}=\ldots=\a_m=0.
   $$
Remarquons que les $\psi_n$ sont lin{\'e}airement ind{\'e}pendants, car
la projection de $\psi_n$ sur $x_1^{-p}y_1^{-p}$-relativement
aux autres {\'e}l{\'e}ments de la base 
$\displaystyle\prod\limits_{-\infty}^{+\infty}x_i^{\a_i}y_i^{\b_i}$ 
o{\`u} $(\a_i)$ et $(\b_i)$ sont deux suites presque nulles dans 
${\ZZ}^{\ZZ}$-vaut $q^{{n\over n+1}}\delta_{n}^{p}$ o{\`u}
$\delta_{n}^{p}$ d{\'e}signe le symbole de Kronecker.

Notons {\'e}galement que pour tout $n,\, \deg{\psi_n}=0$.
Nous pouvons maintenant {\'e}noncer le th{\'e}or{\`e}me suivant~:
   
  \begin{thm}\label{thprin}
   Les $\psi_n$ sont des densit{\'e}s d'int{\'e}grales de mouvement 
   pour le syst{\`e}me
   de sinus-Gordon quantique sur r{\'e}seau. On a~:
   \begin{align}
    \label{thprina} V_{+}.\psi_n&=0,\\
    \label{thprinb} V_{-}.\psi_n&=-(T-Id)(V_{-}.A_n).
   \end{align}  
  \end{thm}

Les $\psi_n$ peuvent {\^e}tre regroup{\'e}es dans une s{\'e}rie 
g{\'e}n{\'e}ratrice de la fa\c con suivante~:
l'alg{\`e}bre $A_q$ {\'e}tant un $\CC[q,q^{-1}]$-module libre,
on a 
  $$
   A_q\subset \CC(q)\otimes_{\CC[q,q^{-1}]}A_q.
  $$
Soit $S$ l'anneau des s{\'e}ries en $\lambda^{-1}$ {\`a}
coefficients 
dans $\CC(q)\otimes_{\CC[q,q^{-1}]}A_q$. Il existe 
une valuation sur $S$
telle que $v(\lambda^{-1})=1$ car $A_q$ est int{\`e}gre.
L'alg{\`e}bre $S$ est le compl{\'e}t{\'e} de
$\bigl( \CC(q)\otimes_{\CC[q,q^{-1}]}A_q\bigr)[\lambda^{-1}]$
pour $v$. Notons $S_1$ le sous-anneau de $S$ constitu{\'e} des s{\'e}ries
formelles $f$ telles que $v(f-1)\geq 1$. Tout {\'e}l{\'e}ment de $S_1$
poss{\`e}de un inverse qui est $\sum_k (1-f)^k$. Consid{\'e}rons la fonction 
${\ln}_q$~:
   \begin{equation}\label{log}
    \begin{array}{rcl}
     {\ln}_q:\quad S_1&\longrightarrow&S\\
                     f&\longmapsto&
                     \displaystyle\sum\limits_{p=1}^{+\infty}
                     \cfrac{1}{[p]}(1-f^{-1})^p.
    \end{array}
   \end{equation}
Cette s{\'e}rie existe car si $f\in S_1$, alors $f^{-1}\in S_1$.
Dans le cas o{\`u} $q\rightarrow 1$, on retrouve le $\ln$ classique.
On a le th{\'e}or{\`e}me suivant~:

  \begin{thm}\label{secthm}
   Soient~:
   \begin{multline}
      U:=\Biggl(1-(x_1 y_1)^{-1}\lambda^{-1}\biggl(
      1-(y_1 x_2)^{-1}\lambda^{-1}\Bigl(\ldots\\
      \quad\ldots 1-(x_{N-1} y_{N-1})^{-1}\lambda^{-1}\bigl(
      1-(y_{N-1} x_{N})^{-1}\lambda^{-1}\bigr)^{-1}\ldots
      \Bigr)^{-1}\biggr)^{-1}\Biggr)^{-1},
    \end{multline}
et
   \begin{multline}
      V:=\Biggl(1-(y_1 x_2)^{-1}\lambda^{-1}\biggl(
      1-(x_2 y_2)^{-1}\lambda^{-1}\Bigl(\ldots\\
      \ldots 1-(y_{N-1} x_{N})^{-1}\lambda^{-1}\bigl(
      1-(x_{N} y_{N})^{-1}\lambda^{-1}\bigr)^{-1}\ldots
      \Bigr)^{-1}\biggr)^{-1}\Biggr)^{-1}.
   \end{multline}
Alors $U$ et $V$ sont bien d{\'e}finies et appartiennent {\`a}
$S_1$.
De plus, on a~:
   \begin{equation}
    {\ln}_q U+{\ln}_q V=\displaystyle\sum\limits_{p=1}^{+\infty}
    \cfrac{\psi_p}{[p]} \lambda^{-p}\,\pmod{\lambda^{-(2N-1)}}
   \end{equation}  
 \end{thm}

La d{\'e}monstration des th{\'e}or{\`e}mes \ref{thprin} et \ref{secthm} 
fera l'objet de la section {\bf 5}.

Regardons maintenant la traduction de ces th{\'e}or{\`e}mes
au niveau classique.

\section{Le cas classique}

\noindent
Les $\psi_n$ ont bien s{\^u}r des analogues classiques.
Posons~:
\begin{align}
     \psi'_{n,\textnormal{cl}}&:=A'_{n,\textnormal{cl}}
     +B'_{n,\textnormal{cl}},\\
     \intertext{avec:}
     B'_{n,\textnormal{cl}}&:=T^{{1\over 2}}A'_{n,\textnormal{cl}},
   \end{align}
et~:
    \begin{equation}
      A'_{n,\textnormal{cl}}
     :=\displaystyle\sum\limits_{{\underline\a}\in I_n}
      \displaystyle{1\over \alpha_1}
      \displaystyle\prod\limits_{i=1}^{2N-3}
      \cnp{\a_{i}+\a_{i+1}-1}{\a_{i+1}}
      (x_1 y_1)^{-\alpha_1}\ldots (y_{N-1}x_N)^{-\alpha_{2N-2}},
    \end{equation}
o{\`u} $N$ est un entier quelconque v{\'e}rifiant~: 
    \begin{equation}
      n\leq 2(N-1),
    \end{equation}
et $I_n$ est l'ensemble des multi-indices 
    \begin{equation}
     {\underline\a}=(\alpha_1,\ldots,\alpha_{2N-2})\in{\NN}^{*}\times
     {\NN}^{2N-3}
    \end{equation}
tels que~:
    \begin{equation}
      \alpha_1+\ldots+\alpha_{2N-2}=n.
    \end{equation}

On a le corollaire du th{\'e}or{\`e}me pr{\'e}c{\'e}dent~:

  \begin{cor}
   Les $\psi'_{n,\textnormal{cl}}$ sont des densit{\'e}s d'
   int{\'e}grales de mouvement
   classiques pour le syst{\`e}me de sinus-Gordon classique sur
   r{\'e}seau. On a~:
     \begin{align}
       \{ \Sigma^+,\psi'_{n,\textnormal{cl}}\}&=0,\\
       \{ \Sigma^-,\psi'_{n,\textnormal{cl}}\}&=
        (T-Id)(\{ A_{n,\textnormal{cl}},\psi'_{n,\textnormal{cl}}\})
     \end{align}
  \end{cor}

On peut rassembler les $\psi'_{n,\textnormal{cl}}$ dans une s{\'e}rie
g{\'e}n{\'e}ratrice~:
   \begin{thm}
    Soit $\lambda^{-1}$ une ind{\'e}termin{\'e}e. Dans le cas classique, on 
    a la formule suivante~:
     $$
      \begin{array}{l}
      \ln\cfrac{1}{1-\cfrac{{(x_1 y_1)}^{-1}\lambda^{-1}}
       {1-\cfrac{{(y_1 x_2)}^{-1}\lambda^{-1}}
       {\ddots\cfrac{\ddots}
       {1-\cfrac{{(x_{N-1} y_{N-1})}^{-1}\lambda^{-1}}
       {1-{(y_{N-1} x_{N})}^{-1}\lambda^{-1}}}}}}
      +\ln\cfrac{1}{1-\cfrac{{(y_1 x_2)}^{-1}\lambda^{-1}}
       {1-\cfrac{{(x_2 y_2)}^{-1}\lambda^{-1}}
       {\ddots\cfrac{\ddots}
       {1-\cfrac{{(y_{N-1} x_{N})}^{-1}\lambda^{-1}}
       {1-{(x_{N} y_{N})}^{-1}\lambda^{-1}}}}}}\\
       =\displaystyle\sum\limits_{n=1}^{2N-2}
       \psi'_{n,\textnormal{cl}}
       \lambda^{-n}\,\pmod{\lambda^{-(2N-1)}}.
      \end{array}
     $$
   \end{thm}

Le fait que les coefficients du d{\'e}veloppement
du logarithme des fractions continues du th{\'e}or{\`e}me 
pr{\'e}c{\'e}dent donne des densit{\'e}s d'int{\'e}grales de mouvement 
du syt{\`e}me de sinus-Gordon classique sur r{\'e}seau a {\'e}t{\'e} d{\'e}montr{\'e} par 
B. Enriquez et B. Feigin par des m{\'e}thodes cohomologiques
(\cite{ENR}).

Revenons au cas quantique, et {\'e}valuons $q$.

\section{Sp{\'e}cialisation de $q$ en un nombre complexe distinct de 1}

\noindent
En vertu du lemme \ref{lemint}, on peut sp{\'e}cialiser $q$ en un
nombre complexe $q_0$ distinct de $1$.
Soit $\pi$ le morphisme d'anneaux~:
   $$
    \begin{array}{cccc}
      \pi :\quad&A_q&
      \longrightarrow &{\CC}[x_i^{\pm 1},y_i^{\pm 1},\, i\in\ZZ]\\
      &q^{\pm 1}&\longmapsto&q_0^{\pm 1}\\
      &x_i^{\pm 1}&\longmapsto&x_i^{\pm 1}\\
      &y_i^{\pm 1}&\longmapsto&y_i^{\pm 1}
    \end{array}
   $$
Les $\pi(\psi_n)$ conviennent comme densit{\'e}s d'int{\'e}grales de
mouvement classiques, et forment une famille libre.
Un cas int{\'e}ressant {\`a} consid{\'e}rer est celui o{\`u}
$q_0$ est une racine primitive  $l^{{\text{{\`e}me}}}$ de l'unit{\'e}.
Dans ce cas, $X_i=x_i^{l}$ et $Y_i=y_i^{l}$ appartiennent
clairement au centre de $A_{q}$. Il est int{\'e}ressant de remarquer que
les $\pi\bigl(\psi_k(x_i,y_i)\bigr)$  sont proportionnelles aux
int{\'e}grales de mouvement classiques 
$\psi_{k,\textnormal{cl}}(X_i,Y_i)$ 
comme nous allons le voir ci-dessous.
\subsection{Sp{\'e}cialisation de $q$ en une racine de l'unit{\'e}}

\noindent
Soit $q_0\not= 1$ une racine primitive $l^{{\text{{\`e}me}}}$ de $1$.
Notons par $v$ la valuation en $q-q_0$ sur 
${\CC}[q]$, ainsi que sur son corps des fractions ${\CC}(q)$. 
Posons~: 
   $$
    {\CC}(q)_{+}:=\lbrace f\in
    {\CC}(q),\, v(f)\geq 0\rbrace.
   $$
L'application $\pi$ se prolonge {\`a} ${\CC}(q)_{+}$~:
   $$
    \begin{array}{rcl}
     \pi:\quad {\CC}(q)_{+}&\longrightarrow&\CC\\
                          f&\longmapsto&f(q_0) 
     \end{array}
   $$
Le th{\'e}or{\`e}me que nous avons en vue est le suivant~:

  \begin{thm}\label{vue}
    Soit $q_0$ une racine primitive $l^{{\text{{\`e}me}}}$ de l'unit{\'e},
    et $n$ un entier divisible par $l$. Posons $n=ln'$, et
    consid{\'e}rons
    pour $i\in\lbrace 1,\ldots,n\rbrace$, les variables commutatives
    $X_i=x_i^{l}$ et $Y_i=y_i^{l}$.
    Alors, on a~:
     $$
      \pi\bigl(\psi_n(x_i,y_i)\bigr)=
      n'\psi_{n',\textnormal{cl}}(X_i,Y_i).
     $$
   \end{thm}

La d{\'e}monstration de ce th{\'e}or{\`e}me est une simple cons{\'e}quence de
la s{\'e}rie de lemmes suivants~:

\begin{lem}\label{L1}
Soient $a$ et $b$ deux entiers avec $a\equiv b\,\pmod{l}$.
Alors, $$v([a])=v([b])=\left\lbrace
\begin{array}{ll}
0&{\text{si }} l\nmid a;\\
1&{\text{si }} l\mid a.
\end{array}\right.$$
\end{lem}

 \begin{dem}
   Tout vient du fait que le polyn{\^o}me $q^n-1$ est s{\'e}parable.
 \end{dem}

 \begin{lem}\label{L2}
   Sous les m{\^e}mes hypoth{\`e}ses que le lemme \ref{L1},
   on a~: 
     $$
       \pi([a])=\pi([b]).
     $$
   De plus, $\pi([a])=0$ si et seulement si $l\mid a$.
 \end{lem}
\begin{dem}
On a $[a+l]=q^l[a]+[l]$, et le r{\'e}sultat s'ensuit.
\end{dem}
\begin{lem}\label{L3}
Soient $a$ et $b$ deux entiers divisibles par $l$, avec $a=la'$ et
$b=lb'$. Alors, $v\Bigl ( \displaystyle{[a]\over [b]}\Bigr )=0$ et 
$\pi\Bigl ( \displaystyle{[a]\over [b]}\Bigr )
=\displaystyle{a'\over b'}=\displaystyle{a\over b}$.
\end{lem}
\begin{dem}
Le fait que $v\Bigl ( \displaystyle{[a]\over [b]}\Bigr )=0$ provient du
lemme \ref{L1}. Pour le reste, on a 
$$\begin{array}{rcl}
\displaystyle{[a]\over [b]}&=&\displaystyle{q^{a}-1\over q^{b}-1}\\
&=&\displaystyle{{(q^l)}^{a'}-1\over {(q^l)}^{b'}-1}\\
&=&\displaystyle{1+q^l+\ldots+q^{l(a'-1)}\over 
1+q^l+\ldots+q^{l(b'-1)}}
\end{array}$$
D'o{\`u} $\pi\Bigl ( \displaystyle{[a]\over [b]}\Bigr )
=\displaystyle{a'\over b'}$.
\end{dem}
\begin{lem}\label{L4}
Soient $a_1$ et $a_2$ deux entiers, avec $l\mid a_1$.
Alors,
$$v\Bigl (\C{a_1+a_2-1}{a_2}\Bigr )=
\left\lbrace\begin{array}{ll}
0&{\text{si }}l\mid a_2;\\
1&{\text{si }}l\nmid a_2.
\end{array}\right.
$$
\end{lem}
\begin{dem}
On a :
$$\begin{array}{rcl}
\C{a_1+a_2-1}{a_2}&=&\displaystyle{[a_1][a_1+1]\ldots [a_1+a_2-1]
\over [1]\ldots [a_2]}\\
&=&\displaystyle{[a_1]\over [a_2]}\times\displaystyle{
[a_1+1]\ldots [a_1+a_2-1]\over [1]\ldots [a_2-1]}.
\end{array}$$
Donc, en supposant $l\mid a_1$, et en appliquant le lemme \ref{L1},
on obtient~:
$$\begin{array}{rcl}
v\Bigl (\C{a_1+a_2-1}{a_2}\Bigr )&=&v([a_1])-v([a_2])\\
&=&1-v([a_2]).
\end{array}$$
On r{\'e}applique ensuite le lemme \ref{L1}, pour obtenir le r{\'e}sultat.
\end{dem}
\begin{lem}\label{L5}
Sous les m{\^e}mes hypoth{\`e}ses que pour le lemme \ref{L4},
on a~: 
$$\pi\Bigl (\C{a_1+a_2-1}{a_2}\Bigr )=0{\text{ si }}l\nmid a_2.$$
Si $l\mid a_2$, en posant $a_1=l a'_1$ et $a_2=l a'_2$,
on a~:
$$\pi\Bigl ( \C{a_1+a_2-1}{a_2}\Bigr )=\cnp {a'_1+a'_2-1}{a'_2}.$$ 
\end{lem}
\begin{dem}
D'apr{\`e}s le lemme \ref{L4},
il suffit de consid{\'e}rer le cas o{\`u} $a_1$ et $a_2$ sont 
divisibles par $l$,
avec $a_1=l a'_1$ et $a_2=l a'_2$.
On a~:
   $$
    \C{a_1+a_2-1}{a_2}=\displaystyle{[a_1]\over [a_2]}\times
    \displaystyle\prod_{k=1}^{a_2-1}\cfrac{[a_1+k]}{[k]}
   $$
Si $l\nmid k$, alors d'apr{\`e}s le lemme \ref{L2}, on a~:
 $$
  \pi\bigl(\cfrac{[a_1+k]}{[k]}\bigr)=1. 
 $$
De plus, si $k\in\lbrace 1,\ldots,a'_2-1\rbrace$,
alors, d'apr{\`e}s le lemme \ref{L3},
 $$
  \pi\bigl(\cfrac{[a_1+kl]}{[kl]}\bigr)=
  \displaystyle\cfrac{a'_1+k}{k}.
 $$
Le r{\'e}sultat en d{\'e}coule.

\end{dem}
\begin{lem}\label{L6}
Pour tout $x\in\ZZ$, on pose $x'=\displaystyle{x\over l}$ 
si $l\mid x$. Soient $p,a_1,\ldots,a_p$ des entiers tels que
$a_1+\ldots +a_p =n$. alors,
$$\pi\Bigl ( \displaystyle{[n]\over [a_1]}F_q(a_1,\ldots, a_p)\Bigr )
=\left\lbrace
\begin{array}{ll}
0,&{\text{ s'il existe $i$ tel que }}l\nmid a_i;\\
\displaystyle{n'\over a'_1}F(a'_1,\ldots,a'_p),&{\text{ sinon}.}
\end{array}\right.$$
avec comme pr{\'e}c{\'e}demment, 
$F(a'_1,\ldots,a'_p)=\displaystyle
\prod\limits_{i=1}^{p-1}\cnp{a'_1+a'_2-1}{a'_2}$.
\end{lem}
\begin{dem}
Notons
$A:=\displaystyle{[n]\over [a_1]}F_q(a_1,\ldots, a_p).$
Le fait que $A$ ait une image par $\pi$ est une cons{\'e}quence du lemme
\ref{lemint}. Le lemme d{\'e}coule des {\'e}tapes
{\'e}l{\'e}mentaires suivantes~:

\begin{enumerate}

\item 
Si $l\nmid a_1$, alors d'apr{\`e}s le lemme \ref{L1}, $v(A)\geq 1$
car $l\mid n$ et $F_q(a_1,\ldots, a_p)\in\ZZ [q]$.
Donc, $\pi (A)=0$.\\

\item
Si $l\mid a_1$ et $l\nmid a_2$, on {\'e}crit 
  $$
   A=\displaystyle{[n]\over [a_1]}
     \C{a_1+a_2-1}{a_2}F_q(a_2,\ldots,a_p),
  $$
et $v(A)\geq 1$ d'apr{\`e}s le lemme \ref{L4}.
Donc, $\pi (A)=0$.\\

\item 
Par suite, de proche en proche, on montre que $v(A)\geq 1$ (et donc
$\pi (A)=0$) s'il existe $i$ tel que $l\nmid a_i$.\\

\item
Si tous les $a_i$ sont divisibles par $l$, le r{\'e}sultat s'obtient 
imm{\'e}diatement {\`a} partir du lemme \ref{L3} et du lemme \ref{L5},
apr{\`e}s avoir {\'e}crit~:
   $$
    A=\displaystyle{[n]\over [a_1]}
    \prod\limits_{i=1}^{p-1}\C{a_i+a_{i+1}-1}{a_{i+1}}.
   $$
\end{enumerate}
\end{dem}

Par suite, nous avons prouv{\'e} le th{\'e}or{\`e}me \ref{vue}.

\section{D{\'e}monstration des th{\'e}or{\`e}mes \ref{thprin} et
\ref{secthm}}
\noindent
Nous aurons besoin des formules suivantes~:
\begin{align}
\label{3.11} (x_r y_r)^{-n}&=q^{\lbrace n\rbrace }x_r^{-n} y_r^{-n}\\
\label{3.12} (y_{r-1} x_r)^{-n}&=q^{\lbrace n\rbrace }y_{r-1}^{-n} x_r^{-n}\\
\label{3.15} x_r^p (x_r y_r)^n &= q^{-np} (x_r y_r)^n x_r^p\\
\label{3.19} x_r^p (y_{r-1} x_r)^n &= q^{np} (y_{r-1} x_r)^n x_r^p
\intertext{avec par convention,}
    \lbrace n\rbrace&=\displaystyle\sum\limits_{i=1}^{n}i=
    \cfrac{n(n+1)}{2}.
\end{align}
\subsection{D{\'e}monstration du th{\'e}or{\`e}me \ref{secthm}~:}

\noindent
Soit $B_q^{(n)}$ l'alg{\`e}bre engendr{\'e}e sur $\CC(q)$ par les
g{\'e}n{\'e}rateurs $t_1,\ldots,t_n$ et les relations~:
   $$
    (R_n):\begin{cases}
         t_{i+1} t_i=q t_i t_{i+1}&\textnormal{pour }i\in\{
         1,\ldots,n-1\} ;\\
         t_i t_j =t_j t_i &\textnormal{pour }\mid i-j\mid \geq 2.
         \end{cases}
   $$
On montre que $B_q^{(n)}$ est int{\`e}gre, et qu'il existe une
valuation $v$ sur $B_q^{(n)}$ en posant $v(t_i)=1\, \forall i$.
On note $T$ le compl{\'e}t{\'e} de $B_q^{(2N-2)}$ pour la topologie
ultram{\'e}trique d{\'e}finie par $v$. C'est le quotient de l'alg{\`e}bre
non commutative des s{\'e}ries formelles en $t_1,\ldots,t_n$
par l'id{\'e}al engendr{\'e} par les relations $(R_{2N-2})$. On d{\'e}finit 
{\'e}galement $T_1$ comme {\'e}tant le sous-anneau de $T$ form{\'e} des
$f\in T$ tels que $v(f-1)\geq 1$. D'apr{\`e}s les r{\`e}gles de
$q-$commutation, on dispose
de deux morphismes naturels~:
   $$
    \begin{array}{rcl}
    \phi:\quad B_q^{(2N-2)}&\longrightarrow&
        \bigl( \CC(q)\otimes_{\CC[q,q^{-1}]}A_q\bigr) [\lambda^{-1}]\\
    \textnormal{si }1\leq k\leq N-1,\quad
        t_{2k-1}&\longmapsto&(x_k y_k)^{-1}\lambda^{-1}\\
        t_{2k}&\longmapsto&(y_k x_{k+1})^{-1}\lambda^{-1}
    \end{array}
   $$
et $\psi$~:
   $$
    \begin{array}{rcl}
      \psi:\quad B_q^{(2N-2)}&\longrightarrow&
    \bigl( \CC(q)\otimes_{\CC[q,q^{-1}]}A_q\bigr) [\lambda^{-1}]\\
  \textnormal{si }1\leq k\leq N-1,\quad
     t_{2k-1}&\longmapsto&(y_k x_{k+1})^{-1}\lambda^{-1}\\
     t_{2k}&\longmapsto&(x_{k+1} y_{k+1})^{-1}\lambda^{-1}
    \end{array}
   $$
Il est clair que les applications $\phi$ et $\psi$ sont continues pour
$v$. Par suite, on peut prolonger $\phi$ et $\psi$ en 
$\widehat{\phi}$ et $\widehat{\psi}$  de $T$ dans $S$, et 
$T_1$ est envoy{\'e} dans $S_1$ par $\widehat{\phi}$ 
ou $\widehat{\psi}$.

\begin{lem}\label{uno}
Soit 
   $$
    M_{n}(t_1,\ldots,t_n):=
    \Biggl(1-t_1\biggl(1-t_2\Bigl(1-
    \ldots\bigl(1-t_{n-1}(1-t_n)^{-1}\bigr)^{-1}
    \ldots\Bigr)^{-1}\biggr)^{-1}\Biggr)^{-1}.
   $$
Alors $M_{n}(t_1,\ldots,t_n)$ existe et 
$M_{n}(t_1,\ldots,t_n)\in T_1$
\end{lem}
\begin{dem}
Par r{\'e}currence sur $n$.

\noindent
$\bullet$ Si $n=1$, le r{\'e}sultat est clair.

\noindent
$\bullet$ Si $n\geq 2$, les variables $t_2,\ldots,t_n$ satisfont les
m{\^e}mes relations que les varialbles $t_1,\ldots,t_{n-1}$ 
(et r{\'e}ciproquement). Donc, par hypoth{\`e}se de r{\'e}currence,
$v\bigl(M_{n-1}(t_2,\ldots,t_n)\bigr)\geq 0.$
Donc, $1-t_1M_{n-1}(t_2,\ldots,t_n)\in T_1$,
et,
   $$
    M_n(t_1,\ldots,t_n)=
    \bigl(1-t_1M_{n-1}(t_2,\ldots,t_n)\bigr)^{-1}
   $$
existe et appartient {\`a} $T_1$.
\end{dem}

\begin{lem}\label{dos}
Soit $(u,v)\in T\times T_1$ tel que $vu=quv$. Alors,
   $$
    \forall n\in\NN,\quad \bigl[ u(1-v)^{-1}\bigr]^{n}=
    \displaystyle\sum\limits_{k\geq 0}\C{n+k-1}{k}u^{n}v^{k}.
   $$
\end{lem}
\begin{dem}
Une simple r{\'e}currence nous montre que 
   $$
    \forall n\in\NN,\quad \bigl[ u(1-v)^{-1}\bigr]^{n}=
    u^{n}(1-v)(1-qv)\ldots (1-q^{n-1}v).
   $$
Le r{\'e}sultat d{\'e}coule alors de l'identit{\'e} classique~:
   $$
    \displaystyle\prod\limits_{s=0}^{n-1}(1-q^{s}X)=
     \displaystyle\sum\limits_{k\geq 0}\C{n+k-1}{k}X^k.
   $$
\end{dem}

Sur $T_1$, on peut d{\'e}finir une fonction ${\ln}_q$ par la m{\^e}me formule
que (\ref{log}). On a le r{\'e}sultat suivant~:
\begin{lem}\label{tres}
Soit $n\in{\NN}^*$. On a~:
   $$
    {\ln}_q\bigl[ M_n(t_1,\ldots,t_n)\bigr]=
    \displaystyle\sum\limits_{\a_1>0,\a_2,\ldots,\a_n}
    \cfrac{1}{[\a_1]}F_q(\a_1,\ldots,\a_n)t_1^{\a_1}\ldots t_n^{\a_n}.
   $$
\end{lem}
\begin{dem}
L'expression a bien un sens d'apr{\`e}s le lemme \ref{uno}.
Remarquons que la formule d{\'e}finissant la fonction ${\ln}_q$ permet
d'{\'e}crire~:
   $$
    \forall f\in T_1,\quad {\ln}_q\bigl[{(1-f)}^{-1}\bigr]=
    \displaystyle\sum\limits_{k\geq 1}\cfrac{1}{[k]} f^k.
   $$
On a 
   $$
    M_n(t_1,\ldots,t_n)=\bigl( 1-t_1 M_{n-1}(t_2,\ldots,t_n)\bigr)^{-1}.
   $$
Donc, 
   $$
    {\ln}_q\bigl[M_n(t_1,\ldots,t_n)\bigr]=
    \displaystyle\sum\limits_{\a_1>0}\cfrac{1}{[\a_1]}
    \bigl(t_1 M_{n-1}(t_2,\ldots,t_n)\bigr)^{\a_1}.
   $$
Or,
   $$
    M_{n-1}(t_2,\ldots,t_n)=
    \bigl(1-t_2 M_{n-2}(t_3,\ldots,t_n)\bigr)^{-1},
   $$
et $t_2$ commute avec $M_{n-2}(t_3,\ldots,t_n)$. Donc, d'apr{\`e}s le
lemme \ref{dos}, on a~:
   $$
     {\ln}_q\bigl[M_n(t_1,\ldots,t_n)\bigr]=
     \displaystyle\sum\limits_{\a_1>0,\a_2\geq 0}\cfrac{1}{[\a_1]}
     \C{\a_1+\a_2 -1}{\a_2}t_1^{\a_1}
     \bigl(t_2 M_{n-2}(t_3,\ldots,t_n)\bigr)^{\a_2}.
   $$
Le r{\'e}sultat s'obtient alors en utilisant de mani{\`e}re r{\'e}p{\'e}t{\'e}e
le lemme \ref{dos}.
\end{dem}

On ach{\`e}ve ensuite la d{\'e}monstration du th{\'e}or{\`e}me \ref{secthm}
en remarquant que les fonctions ${\ln}_q$ sur $T_1$ et $S_1$ commutent
avec $\widehat{\phi}$ et $\widehat{\psi}$, et que
   \begin{align*}
               \widehat{\phi}\bigl(M_n(t_1,\ldots,t_n)\bigr)&=U,\\
\text{et }\quad\widehat{\psi}\bigl(M_n(t_1,\ldots,t_n)\bigr)&=V.
   \end{align*}

\subsection{D{\'e}monstration du th{\'e}or{\`e}me \ref{thprin}}

\noindent
L'{\'e}galit{\'e} (\ref{thprinb}) d{\'e}coule de (\ref{thprina}). En effet,
si $V_+ .\psi_n=0$, alors, avec les notations de la section {\bf 2},
   $$
    [\Sigma^+,A_n]+[\Sigma^+,T^{{1\over 2}}A_n]=0,
   $$
Or, l'application $T^{{1\over 2}}$ est un isomorphisme
d'alg{\`e}bres qui transforme $\SP$ en $\SM$ (ou plut{\^o}t,
qui permute $V_+$ avec $V_-$). Donc,
   $$
    [\Sigma^{-},T^{{1\over 2}}A_n]+[\Sigma^{-},T A_n]=0.
   $$
Donc, 
   $$
    [\Sigma^-,B_n]=-T[\Sigma^-,A_n],
   $$
car $T\Sigma^-=\Sigma^-$ (ou $T\circ V_- =V_-\circ T$).
Donc, 
   \begin{align*}
   V_{-}.\psi_n&=[\Sigma^{-},A_n]+[\Sigma^-,B_n]\\
               &=[\Sigma^{-},A_n]-T\bigl( [\Sigma^{-},A_n]\bigr)\\
               &=-(T-Id)(V_- .A_n).
   \end{align*}

\bigskip
\noindent
Avant de commencer la d{\'e}monstration de (\ref{thprina}),
notons que dans la d{\'e}finition de $A_n$, on peut remplacer
l'ensensemble $I_n$ par les ensembles $J_n$ ou $J''_n$, avec~:
   \begin{align}
   J_n&=\left\lbrace (\alpha_1,\ldots ,\alpha_{2N-2})\in 
     {\NN}^{*}\times {\ZZ}^{2N-3}/\; \alpha_1+\ldots +\alpha_{2N-2}=n
      \right\rbrace,\\
   J''_n&=\left\lbrace (\alpha_1,\ldots ,\alpha_{2N-2})\in 
         {\NN}^{*}\times{\NN}\times {\ZZ}^{2N-4}/\; 
         \alpha_1+\ldots +\alpha_{2N-2}=n\right\rbrace.
   \end{align}
Ceci vient du fait suivant qui est une cons{\'e}quence de la
d{\'e}finition de la fonction $F_q$~:

\medskip
\noindent
{\bf Fait 2~:} Soient $m$ et $(\a_1,\ldots,\a_m)\in\ZZ^m$.
S'il existe $i$ tel que $\a_i<0$, alors $F_q(\a_1,\ldots,\a_m)=0$.

\bs
Nous poserons {\'e}galement~:
   $$
    J'_n= \left\lbrace (\alpha_1,\ldots ,\alpha_{2N-2})\in 
    {\NN}\times {\ZZ}^{2N-3}/\; \alpha_1+\ldots +\alpha_{2N-2}=n
    \right\rbrace.
   $$
Ecrivons~:
   \begin{equation}
    \label{4.2} [\Sigma^{+},\psi_n]=\left[x_1,A_n\right]+
    \bigl[\displaystyle\sum\limits_{k\geq 2}
    x_k,A_n\bigr]+\left[\Sigma^{+},B_n\right],                      
   \end{equation}
et calculons s{\'e}par{\'e}ment chacun de ces termes.

\bs
\noindent
{\bf 1. Calcul de $\left[x_1,A_n\right]$.}

\bs
\noindent En utilisant~(\ref{3.15}), il vient~:
\begin{multline*}
\left[x_1,A_n\right]\\
=\SIA\dfrac{1-q^{-\ab}}{[\ab]}\,F_q(\ab,\ldots,\af)\:x_1\,\HT
\end{multline*}
Donc, en s{\'e}parant les cas $\ab =1$ et $\ab >1$ et 
{\`a} l'aide de l'{\'e}galit{\'e}~:
\begin{equation}\label{4.1}
\forall n>1,\; 
\dfrac{1-q^{-n}}{\q{n}}=\dfrac{1}{\q{n-1}}(q^{-1}-q^{-n}),
\end{equation} 
on obtient
\begin{equation}\label{4.3}
\left[ x_1,A_n\right]=R_1+R_2
\end{equation}
avec
\begin{multline*}
R_1=(1-q^{-1})\\
\times\SJ\,F_q(1,\,\ac,\ldots,\af)x_1\HV
\end{multline*}
et
\begin{multline}\label{4.4}
R_2=\SK\dfrac{1}{[\ab -1]}(q^{-1}-q^{-\ab})\, F_q(\ab,\ldots,\af)\\
\times x_1\HT
\end{multline}
Remarquons que $R_1$ peut encore s'{\'e}crire
\begin{multline}\label{4.5}
R_1=\SID\dfrac{q}{\q{\ac +1}}(q^{-1}-q^{-\a_1})\, F_q(\ab,\ldots,\af)\\
\times x_1\HV
\end{multline}
\bs
{\bf 2. Calcul de $\bigl[\displaystyle\sum\limits_{k\geq
2}x_k,A_n\bigr]$.}

\bs
\noindent
En utilisant~(\ref{3.15}) et (\ref{3.19}),
et, en convenant de poser que ${\a}_{2N-1}=0$ pour 
${\underline\a}\in J_n$, on a
\begin{multline*}
\bigl[\displaystyle\sum\limits_{k\ge 2} x_k,A_n\bigr]=
\displaystyle\sum\limits_{k=2}^{N}\,\SIA\dfrac{1}{[\ab]}
(q^{-\ai}-q^{-\aj})F_q(\ab,\ldots,\af)\\
\times \hb\ldots \hi\,x_k\,\hj\ldots\hf
\end{multline*}
Mais, d'apr{\`e}s~(\ref{3.11}) et (\ref{3.12}), on a 
\begin{multline*}
\hb\hc\ldots \hi\,x_k=q^{X_k}\\
 \times\, x_1{\ob}^{-(\ab +1)}{\oc}^{-(\ac -1)}\ldots {\oi}^{-(\ai -1)}
\end{multline*}
avec
\begin{equation*}
\begin{split}
X_k&=\lbrace \ab\rbrace +\cdots+\lbrace\ai\rbrace 
-\left( \lbrace \ab+1\rbrace +\lbrace\ac-1\rbrace
+\cdots+\lbrace\ai-1\rbrace\right)\\
&=-(\ab +1)+\ac-\ldots +\ai
\end{split}
\end{equation*}
En outre, on note que pour tout $k\in\lbrace 2,\ldots,N\rbrace$,
l'application
$$
f_k:\qquad\a_i\longrightarrow
\left\lbrace
\begin{array}{cl}
\a_i +1&\textnormal{si}\;i\; \textnormal{impair }\le 2k-3\\
\a_i -1&\textnormal{si}\;i\; \textnormal{pair }\le 2k-2\\ 
\a_i&\textnormal{sinon}
\end{array}
\right.
$$
est une bijection de $J_n$ sur 
$\left\lbrace {\underline{\a}}\in J_n;\; 
a_{1} >1\,\right\rbrace$. Il en r{\'e}sulte que
\begin{multline}\label{4.6}
[\displaystyle\sum\limits_{k\ge 2} x_k,A_n]=\SK\displaystyle
\sum\limits_{k=2}^{N}
\dfrac{q^{Y_k}}{\q{\ab -1}}\left( q^{-(\ai +1)}-q^{-\aj}\right)\\
\times\, F_q(\ab -1,\,\ac +1,\cdots,\ai +1,\,\aj,\cdots,\,\af)\\
\times\, x_1\HT
\end{multline}
avec
\begin{align}
\notag Y_k&=-\ab+(\ac +1)-\cdots +\ai +1\\
\label{4.7} &=k-1+\sum\limits_{j=1}^{2k-2}(-1)^j{\a}_j 
\end{align}
et en gardant la convention $\ag =0$ pour ${\underline{\a}}\in J_n$.

\bs
\noindent
{\bf 3. Calcul de $\left[\SP,B_n\right]$.}

\bs
\noindent
On proc{\`e}de comme dans le cas pr{\'e}c{\'e}dent. 
D'apr{\`e}s le formulaire, $x_k$
et $B_n$ commutent si $k\notin\lbrace 2,\ldots,N\rbrace$. On a donc
\begin{multline*}
[\SP,B_n]=\displaystyle
\sum\limits_{k=2}^{N}\SIB\dfrac{1}{[\bc]}
(q^{-\bi}-q^{-\bj})F_q(\bc,\ldots,\bg)\\
\times\, \gc\ldots \gi\,x_k\,\gj\ldots (x_N y_N)^{-\beta_{2N-1}}
\end{multline*}
 Les relations~(\ref{3.11}) et (\ref{3.12}) montrent que
\begin{multline*}
\gc\ldots \gi\,x_k=q^{X'_k}\\
\times\, x_1{\ob}^{-1}{\oc}^{-(\bc -1)}{(x_2 y_2)}^{-(\b_3 +1)}
\ldots {\oi}^{-(\bi -1)}
\end{multline*}
avec 
\begin{align*}
X'_k&=\lbrace \b_2\rbrace +\cdots +\lbrace \b_{2k-2}\rbrace
-\left( 1 + \lbrace\bc -1\rbrace+\lbrace\b_3 +1\rbrace
+\cdots+\lbrace \bi+1\rbrace\right)\\
&=-1+\b_2-(\b_3+1)-\cdots+\bi.
\end{align*}
De plus, l'application
$$
\begin{array}{ccc}
J_n&\longrightarrow&J'_{n-1}\\
(\bc,\ldots,\bg)&\longrightarrow &(\ac,\ldots,\ag)\\
\b_i&\longrightarrow &
\left\lbrace
\begin{array}{cl}
\b_i +1&\textnormal{si}\;i\; \textnormal{impair }\le 2k-3\\
\b_i -1&\textnormal{si}\;i\; \textnormal{pair }\le 2k-2\\ 
\b_i&\textnormal{sinon}
\end{array}
\right.
\end{array}
$$
est bijective.
Donc, 
\begin{multline}\label{4.8}
[\SP,B_n]=\displaystyle
\SIC\sum\limits_{k=2}^{N}\dfrac{q^{Y'_k}}{\q{\ac +1}}
(q^{-(\ai+1)}-q^{-\aj})\\
\times\, F_q(\ac+1,\a_3-1,\ldots,\ai+1,\aj,\ldots,\ag)\\
\times\, x_1{\ob}^{-1}\hc\ldots (x_N y_N)^{-\alpha_{2N-1}}
\end{multline}
\begin{align*}
{\textnormal{avec }}
{\underline{\a}}=(\ac,\ldots,\ag)
{\textnormal{ et }}Y'_k&=-1+(\ac +1)-\cdots +\ai +1\\
&=k-2+\sum\limits_{j=2}^{2k-2}(-1)^j{\a}_j 
\end{align*}

\medskip
\noindent Arr{\^e}tons nous un instant sur l'expression~(\ref{4.8}),
et, montrons que si un multi-indice $({\underline{\a}};k)$
avec ${\underline{\a}}=(\ab,\ldots,\ag)\in J'_{n-1}$
v{\'e}rifie $\ag\ge 1$, alors son terme correspondant est nul.\\
En effet, si $\af\ge 1$, alors, d'apr{\`e}s les faits 1 et 2,
pour que l'on ait
   $$
    F_q(\ac+1,\a_3-1,\ldots,\ai+1,\aj,\ldots,\ag)\not=0,
   $$
il faut que tous les arguments soient $>0$.\\
\smallskip
Donc, 
   $$
    \left\lbrace
     \begin{array}{c}
      \ac\ge 0\\
      \a_3\ge 2\\
      \vdots\\
      \a_{2k-3}\ge 2\\
      \ai\ge 0\\
      \aj\ge 1\\
      \vdots\\
      \ag\ge 1
     \end{array}
    \right.
   $$
et
   $$
    n-1=\sum\a_i\ge 2(k-2)+2(N-k)+1=2N-3.
   $$\\
\smallskip
Or, $n\le 2(N-1)$. Donc toutes les in{\'e}galit{\'e}s pr{\'e}c{\'e}dentes
sont des {\'e}galit{\'e}s, et, en particulier, on a $\ai =0$
et $\aj=1$. D'o{\`u}, $(q^{-(\ai +1)}-q^{-\aj})=0$ et le terme
correspondant {\`a} $({\underline{\a}};k)$ est nul.
\\
Ainsi, dans~(\ref{4.8}), on peut se restreindre {\`a} sommer sur les
multi-indices $({\underline{\a}};k)$, avec $\ag=0$.
Par suite, on peut {\'e}crire~:
\begin{multline}\label{4.9}
[\SP,B_n]=\displaystyle\SID
\sum\limits_{k=2}^{N}\dfrac{q}{\q{\ac +1}}q^{Y_k}
(q^{-(\ai +1)}-q^{-\aj})\\
\times\, F_q(\ab -1,\ac +1,\cdots,\ai +1,\aj,\cdots,\af)\\
\times\, x_1{\hbab}\hc\cdots\hf
\end{multline}
avec $Y_k$ d{\'e}fini comme dans la formule~(\ref{4.7}) du {\it B}, et,
par convention, 
$\ag =0$ pour ${\underline{\a}}=(\ab,\ldots,\af)
\in J''_{n+1}$.

\bs
\noindent
{\bf 4. Rassemblons nos r{\'e}sultats.}

\bs
\noindent
En regroupant les relations~(\ref{4.2}),
(\ref{4.3}), (\ref{4.4}), (\ref{4.5}), (\ref{4.6}), (\ref{4.9}), 
on obtient
\begin{multline}\label{4.10}
[\SP,P_n]=\SID\dfrac{q}{\q{\ac +1}}\left( \displaystyle\sum\limits_{k=1}^{N}
\Phi_{N-1} (k,\ab,\cdots,\af)\right)\\
\times\,  x_1{\ob}^{-1}\hc\cdots\hf\\
+\SK\dfrac{1}{\q{\ab -1}}\left( \displaystyle\sum\limits_{k=1}^{N}
\Phi_{N-1} (k,\ab,\cdots,\af)\right)\\
\times\, x_1\HT.
\end{multline}
avec $Y_k$ d{\'e}fini comme dans~(\ref{4.7}) et
\begin{multline}
\label{4.11} \Phi_{N-1} (k,\ab,\cdots,\af)=q^{Y_k} 
(q^{-(\ai +1)}-q^{-\aj})\\
\times\, F_q(\ab -1,\ac +1,\cdots,\ai +1,\aj,\cdots,\af)
\end{multline}
et la convention que $\a_0=\ag=0$ pour 
${\underline{\a}}\in J''_{n+1}$ comme pour
${\underline{\a}}\in J_n$.

\bs
\noindent
Soit ${\underline{\a}}=(\ab,\cdots,\af)\in{\ZZ}^{2N-2}$
avec $\a_1\ge 2$, et supposons qu'il existe 
$k\in\,\lbrace 1,\ldots,N\rbrace $ tel que 
$\Phi_{N-1}(k,{\underline{\a}})\not= 0$.
Alors, d'apr{\`e}s la d{\'e}finition de $F_q$, tous les coefficients 
de $f_k({\underline{\a}})$ sont $\ge 0$.
Donc, en particulier, 
$$\ai +1\ge 0\, ,\,\aj\ge 0\, ,\ldots,\,\af\ge 0.$$
On est dans l'un des deux cas suivants~:
\begin{itemize}
\item[1)] $\aj\ge 1$. Alors,
$F_q(\ab -1,\ac +1,\cdots,\ai +1,\aj)\not=
0$ entra{\^\i}ne~ (d'apr{\`e}s les fait 1 et 2):
   $$
    \left\lbrace
     \begin{array}{c}
      \ab-1\ge 1\\
      \ac +1\ge 1\\
      \vdots\\
      \ai+1\ge 1\\
      \aj\ge 1
     \end{array}
    \right.
   $$
ce qui montre que tous les $\a_i$ sont $\ge 0$.
\item[2)] $\aj =0$. Dans ce cas, on a n{\'e}cessairement $\ai +1\not=0$, 
car sinon, $(q^{-(\ai +1)}-q^{-\aj})=0$ et 
$\Phi_{N-1}(k,{\underline{\a}})=0$. Donc $\ai +1\ge 1$. Mais ceci
entra{\^\i}ne comme pr{\'e}c{\'e}demment, ${\a}_{i}\ge 0\quad \forall i$.
\end{itemize}
Par suite, dans les deux sommes de~(\ref{4.10}), on peut se restreindre
aux multi-indices ${\underline{\a}}=(\ab,\cdots,\af)$ tels que
${\underline{\a}}\in{\NN}^{2N-2}-\lbrace 0\rbrace$.
\\
Pour conclure, il suffit donc de prouver le lemme suivant:

\bs
\begin{lem}\label{le1} 
Soient $N\in{\NN}^*$, et $r\in\{ 1,\ldots,N+1\}$. On pose~:
   $$
    \begin{array}{rcl}
     \Phi_{N,r}:{\NN}^{2N}&\longrightarrow&
     \CC[q,q^{-1}]\\
     (a_1,\ldots,a_{2N})&\longmapsto&
     q^{Y_r(a_1,\ldots,a_{2N})}
     \bigl(q^{-(a_{2r-2}+1)}-q^{-a_{2r-1}}\bigr)\\
     &&\times F_q(a_1-1,a_2+1,\ldots,a_{2r-2}+1,a_{2r-1},\ldots,a_{2N})
    \end{array}
   $$
avec la convention que 
   $$
    Y_r(a_1,\ldots,a_{2N})=r-1+\sum\limits_{j=1}^{2r-2}(-1)^j a_j
   $$
et $a_0=a_{2N+1}=0$.
Alors,
   $$
    \forall (a_1,\ldots,a_{2N})\in{\NN}^{2N}-\{ 0\},\quad
    \sum\limits_{k=1}^{N+1}
    \Phi_{N,k}(a_1,\ldots,a_{2N})=0
   $$
\end{lem}
\begin{pr}
La d{\'e}monstration se fait par r{\'e}currence sur $N$.
Elle est laiss{\'e}e au lecteur.
\end{pr}

\noindent
Ceci ach{\`e}ve la d{\'e}monstration du th{\'e}or{\`e}me \ref{thprin}.


{\small\begin{thebibliography}{1}
\bibitem{BABO}
V. Bazhanov, A. Bobenko, N. Reshetikhin,
Comm. Math. Phys., {\bf 175} (1996), 
${\textnormal{n}}^{\circ}2$, 377-400.
\bibitem{BOB}
A. Bobenko, N. Kutz, U. Pinkall,
Phys. Lett., A {\bf 177} (1993), 399-404.
\bibitem{ENR}
B. Enriquez, B. Feigin,
Theor. Math. Phys., {\bf 103} (1995), 738-756.
\bibitem{FAD}
L. Faddeev, A. Volkov,
Theor. Math. Phys., {\bf 92} (1992), 207-214.
\bibitem{FEI}
B. Feigin, E. Frenkel, Lect. Notes in Math., {\bf 1620} (1996),
Springer-Verlag 
\bibitem{IZK}
A.G. Izergin, V.E. Korepin,
Lett. Math. Phys., {\bf 5} (1981), 199-205.
\bibitem{IZKO}
A.G. Izergin, V.E. Korepin,
Nucl.Phys., B {\bf 205} (1985), 401-413.
\bibitem{KAS}
C. Kassel, \emph{Quantum groups}, Springer-Verlag 1994 p.15.
\bibitem{RES}
N. Reshetikhin, Int. Conf. of Math. Phys., 1994, Paris.
\bibitem{SAS}
R. Sasaki, I. Yamanaka, Adv. Stud. in Pure Math., {\bf 16} (1988),
271-296 
\bibitem{ZALZA}
A. Zamolodchikov, Al. Zamolodchikov, Ann. Phys., {\bf 120} (1979),
253-291
\bibitem{ZA} 
A. Zamolodchikov, Theor. Math. Phys., {\bf 65} (1985), 1205
\bibitem{ZAM} 
A. Zamolodchikov, Adv. Stud. in Pure Math., {\bf 19} (1989), 641-674
\end{thebibliography}}
\end{document}